\documentclass[11pt,a4paper]{article}
\usepackage{amsmath,amsfonts,amssymb,color}
\usepackage{enumerate}
\usepackage[english]{babel}

\usepackage[utf8]{inputenc}
\usepackage[T1]{fontenc}

\author{Bruno Luiz Santos Correia and Marc Troyanov}

\date{\today}

\newtheorem{theorem}{\textrm{Theorem}}
\newtheorem{lemma}[theorem]{\rm\bf Lemma}
\newtheorem{corollary}[theorem]{\rm\bf Corollary}

  \newcommand{\qed}{\hfill {\raisebox{-0pt}{\framebox[.23cm]{}}} \newline}

 \newcommand{\rnote}{\textcolor{red}}

\newcommand{\N}{\mathbb{N}}

\newcommand{\Rp}{\mathbb{R}_+}
\DeclareMathOperator{\Card}{Card}

\title{On the Isoperimetric Inequality in Finitely Generated Groups}

\begin{document}

\maketitle

\begin{abstract}\parindent=0pt
We present a sharp version of the isoperimetric inequality for finitely generated groups due to T. Couhlon and L. Saloff-Coste based
on the proof suggested by M. Gromov.  

\medskip 

{\small  {Keywords: Finitely generated groups, growth, isoperimetric inequality, Lambert $W$ function.}}

{\small  {AMS Subject Classification: 20F32;  20F69.}}
\end{abstract}

\section{Introduction and statement of the Main Results}

Let $\Gamma$ be a group generated by the finite symmetric set $S = S^{-1} \subset \Gamma$. We denote by $B_S(n) \subset \Gamma$ the subset of those elements
in $\Gamma$ that can be written as a product of at most $n$ generators in $S$. The  right invariant \textit{word metric} on  $\Gamma$
with respect to the generating set $S$ is   defined as
$$
  d_S(x,y) = \min \{n \in \mathbb{N}  \mid xy^{-1} \in B_S(n)\}.
$$
We shall denote by $\|x\|_S = d_S(e,x)  = \min \{n \in \mathbb{N}   \mid x \in B_S(n)\}$ and  we define the \textit{growth function}   $\gamma_S : \N \to \N$ 
with respect to $S$  as 
$$
  \gamma_S(n) = \Card (B_S(n)).
$$

The \textit{boundary} of a non empty finite subset $D \subset \Gamma$ is defined as
$$
 \partial_S D = \{x\in D \mid \operatorname{dist}_S(x, \Gamma \setminus D) = 1\}
=  \{x\in D \mid \exists s \in S \  \text{ such that } \ sx \not\in D\}.
$$

\medskip

Following T. Couhlon and L. Saloff Coste in \cite{CSC}, we define the ``\textit{reverse growth function}''  of  $(\Gamma, S)$ to be the following  function $\phi_S : \Rp \to \N$:
\begin{equation}\label{def.phi}
 \phi_S(t) = \min \{n \in \N \mid \gamma_S(n) \geq t\}.
\end{equation}

With this preparation, we state  the following isoperimetric inequality:
\begin{theorem} \label{MainTh}
Let $\Gamma$ be a finitely generated group. For any non empty finite subset $D\subset \Gamma$ and any 
real number $\lambda$ such that  $1 < \lambda   \leq  \frac{\Card (\Gamma)}{\Card(D)}$, we have
\begin{equation}\label{mainineq}
  \Card(\partial D)  \geq \left(1 - \frac{1}{\lambda}\right) \frac{\Card(D)}{\phi_S(\lambda \Card( D))}.
\end{equation}
\end{theorem}

\bigskip 

This is the main result of the paper and the proof will be given in the next section. 
It is well known that the growth $\gamma_S(n)$ of an infinite finitely generated group is at least a linear 
and at most an exponential function of $n$. We refer to the survey  \cite{GH} for more on the growth function
of finitely generated groups. 
Optimizing the inequality  (\ref{mainineq}) over $\lambda$ allows us to prove in \S 3 the following explicit isoperimetric inequalities
for groups with a given lower bound on the growth:

\begin{corollary} \label{cortomain}
Let $\Gamma$ be an infinite  finitely generated group with growth function $\gamma_S$.

\smallskip 

\emph{(I)}  If  $\gamma_S(n-1) \geq   Cn^d$ for some constants $C >0$ and  $d \geq 1$ and any integer $n \geq 1$. Then
\begin{equation}\label{isoppol}
  \Card(\partial D)  \geq \frac{C^{\frac{1}{d}} d}{(d+1)^{\frac{d+1}{d}}}\Card(D)^{\frac{d-1}{d}}
\end{equation}
for any finite subset $D\subset \Gamma$.

\medskip 

\emph{(II)} If $\displaystyle \gamma_S(n-1) \geq C{{\rm e}^{b{n}^{\alpha}}}$  for  some constants $C>0$, $b>0$ and
$0 < \alpha \leq 1$,  and  any integer   $n \geq 1$, then 
\begin{equation}\label{isopexp}
  \Card(\partial D)  \geq   \mu(\Card(D))  \cdot   \frac{ \Card(D)}{ \left( \frac {1}{b} (\log (\Card(D))\right)^{1/\alpha}}
 \end{equation}
 for any finite subset $D\subset \Gamma$, 
where $\mu(v)$ is an explicit function such that $\mu(v) \to 1$ as $v \to \infty$.
\end{corollary}

\medskip

\textbf{Remarks.}  
$\circ$ Theorem 1  is a refinement of a result  by  T. Coulhon  and  {L. Saloff-Coste},  see  \cite[th. 1]{CSC}  and   \cite[Theorem 3.2]{PSC}. However the first version of such an isoperimetric inequality, at least for groups of polynomial growth, is usually attributed to N. Varopoulos
(the reader may compare with Property 3 page 388 and its proof page 391 in \cite{Varopoulos}; note that the discussion in this paper is about nilpotent Lie groups.)

\smallskip 

$\circ$ In  Section $\mathrm{E}_+$  of Chapter 6 in the book \cite{Gromov}, M. Gromov stated the inequality (\ref{mainineq}) for $\lambda = 2$ and,  with essentially the same proof, see also \cite{Santos}. We obtain here the general case.  

\smallskip 

$\circ$  The explicit expresion for the function $\mu$ appearing in the isoperimetric inequality \eqref{isopexp} is given in \eqref{defmu2} below.

\section{Proof of Theorem 1.}

Our proof is based on the strategy given by Gromov in \cite{Gromov};  we divide the argument in five steps. 

\smallskip

Step (i).  For any $y \in \Gamma$, we define the following subsets of $D$:
$$
 E_y = \{x \in D \mid yx \not\in D\}  \quad \text{and} \quad I_y = D \setminus E_y = y^{-1}D \cap D.
$$
The points in $E_y$ are the ``exit points'', they are transported outside $D$ by the left translation by $y$, while the points  
in $I_y$ remain inside $D$.   We  claim that 
\begin{equation}\label{SizeEy}
 \Card(E_y) \leq \|y\|_S\Card (\partial D).
\end{equation}
Indeed, suppose that $y = s_n\cdot s_{n-1}\cdots s_1$, with $s_j \in S$ and $n = \|y\|_S$,  and define inductively $y_k \in \Gamma$ by 
$y_0 = e$ and $y_{j} = s_j y_{j-1}$ for $1 \leq k \leq n$.  We then define a map $f : E_y \to \Gamma$
as
\begin{equation} \label{deff}
 f(x) = y_mx, \quad \text{ where } m = \max \{ k \in \N \mid k \leq (n-1) \text{ and } y_k x \in D\}.
\end{equation}
Observe that $f(x) \in \partial D$ for any  $x\in  E_y$ since $yx = y_nx \not \in D$. Furthermore, for any $z \in \partial D$ we have
\begin{equation}\label{preimageZ}
  f^{-1} (z) \subset \{y_0^{-1}z, \dots , y_{n-1}^{-1}z\}.
\end{equation}
 We thus have $\Card (f^{-1}(z)) \leq n$ for any $z \in \partial D$ and therefore 
$$
 \Card(E_y) = \sum_{z\in \partial D} \Card (f^{-1} (z))  \leq n \Card (\partial D) = \|y\|_S \Card (\partial D),
$$
this proves \eqref{SizeEy}.

\medskip

Step (ii). We have 
\begin{equation}\label{SumIy}
  \sum_{y\in \Gamma} \Card (I_y) = \Card(D)^2.
\end{equation}
This follows from the identities\footnote{Permuting the sum in the equality $(*)$ is allowed since  the set of $(x,y) \in \Gamma\times \Gamma$ such that $x\in D$ and $\chi_{{}_D}(yx) \neq 0$ is finite.}
$$
 \sum_{y\in \Gamma} \Card (I_y)  =  \sum_{y\in \Gamma}  \sum_{x\in D}  \chi_{{}_D}(yx)
 \overset{(*)}{=}  \sum_{x\in D}  \sum_{y\in \Gamma}  \chi_{{}_D}(yx)
= \sum_{x\in D}\Card \{y \in \Gamma \mid y \in Dx^{-1}\},
$$
where $\chi_{{}_D}$ is the characteristic function of $D$.

\medskip 

Step (iii). For any $n$ we can find $y \in B_S(n)$ such that 
\begin{equation}\label{SumIy2}
  \gamma_S(n)\Card (I_y)   \leq \Card(D)^2.
\end{equation}
Indeed, choosing  $y \in B_S(n)$ such that   $\Card (I_{y}) = \min \{\Card (I_w)  \mid w \in B_S(n)\}$ and applying 
(\ref{SumIy})  yields
$$
  \gamma_S(n)\Card (I_y) =  \Card(B_S(n))\Card (I_y) \leq \sum_{z\in B_S(n)} \Card (I_z) 
\leq \sum_{z\in \Gamma} \Card (I_z) = \Card(D)^2.
$$

\smallskip 
 
Step (iv). Let us now set  $n = \phi_S(\lambda \Card(D))$, then  $n$ is the smallest integer such that 
$
 \lambda \Card(D)  \leq  \gamma_S(n).
$
From step (iii),  we find $y \in B_S(n)$ such that 
$$
   \lambda \Card(D)\Card (I_{y}) \leq  \gamma_S(n)\Card (I_y) \leq   \Card(D)^2. 
$$
Therefore
\begin{equation}\label{inqyxc}
 \Card (I_{y}) \leq  \frac{1}{\lambda} \Card(D).
\end{equation}

\medskip 

Step (v).  Using \eqref{SizeEy} and \eqref{inqyxc}, we obtain  
\begin{align*}
\left(1-\frac{1}{\lambda}\right) \Card(D) &\leq   \Card(D) - \Card (I_{y}) 
= \Card(E_{y})  \leq \|y\|_S \Card (\partial D)
\\& \leq n  \Card (\partial D)
=  \phi_S(\lambda \Card(D))  \Card (\partial D).
\end{align*}
The proof of Theorem 1 is complete.
\qed

\bigskip

\section{Proof of Corollary \ref{cortomain}}

To prove Corollary \ref{cortomain}, it will be useful to first look at some properties of the function $\phi_S$:
\begin{lemma} The function $\phi_S$ satisfies the following properties
\begin{enumerate}[(i)]
\item $\phi_S$ is a non decreasing function. 
\item $\phi_S$ is a left inverse of $\gamma_S$, that is 
$
  \phi_S \left( \gamma_S (n)\right) = n
$
for any integer $n\in \N$.
\item We have $\gamma_S \left( \phi_S (m)\right) \geq m$ for any $m\in \N$ and 
$
  \gamma_S \left( \phi_S (m)\right) = m
$
for any $m \in \gamma_S(\N)$.
\item Suppose that $ \gamma_S (n) \geq g(n+1)$ for any $n\in \N$, where  $g : [0, \infty)\to [\beta,\infty)$   is a homeomorphism, then 
$\phi_S(t) \leq g^{-1}(t)$ for any $t\in \Rp$.
\end{enumerate}
\end{lemma}

\medskip

\textbf{Remark.}  Simple examples show that the inequality  $ \gamma_S (n) \geq g(n)$ for all $n \in \N$  does \textit{not} generally imply
$\phi_S(m) \leq g^{-1}(m)$  for all $m\in \N$. This is the reason we work with the hypothesis  $ \gamma_S (n) \geq g(n+1)$. 

\medskip

\textbf{Proof.} The  first three statements are elementary, let us prove the last one. 
The  condition  $g (n+1) \leq \gamma_S(n)$  implies 
$$
  \{n \in \N \mid g(n+1) \geq t\} \subset \{n \in \N \mid \gamma_S(n) \geq t\},
$$
therefore
$$
  \phi_S(t) = \min \{n \in \N \mid \gamma_S(n) \geq t \} \leq   \min \{n \in \N \mid g(n+1) \geq t \}.
$$
Because $g(r)$ is monotone increasing, we have $g(n) < g(r) \leq g(n+1)$ whenever $n<r\leq n+1$.
Therefore
$$
   \phi_S(t) \leq   \min \{n \in \N \mid g(n+1) \geq t \} \leq   \min \{r \in \Rp \mid  g(r)\geq t \} =  g^{-1}(t).
$$
\qed

\medskip

Using the previous Lemma and  Theorem \ref{MainTh},  we immediately obtain the following Corollary:

\begin{corollary}\label{corMainIsop}
Let $\Gamma$ be {an infinite}  group generated by a finite symmetric set $S=S^{-1}$. Assume that  $\gamma_S (n) \geq g(n+1)$ for any 
$n\in \N$, where $g : \Rp \to  \Rp$ is a  homeomorphism.  Then the following isoperimetric inequality holds for any non empty finite subset $D \subset \Gamma$:
\begin{equation}\label{ineq2}
  \Card(\partial D)  \geq F(\Card(D)),
\end{equation}
where $F(v)$ is defined as
\begin{equation}\label{defF}
 F(v) = \sup_{1< \lambda < \infty}\left(\frac{\left(1-\tfrac{1}{\lambda}\right)v}{g^{-1}(\lambda v)}\right).
\end{equation}
\end{corollary}
\qed

\bigskip

\textbf{Proof of Corollary \ref{cortomain}.}
We first consider the general situation as in Corollary \ref{corMainIsop}. Let us denote by $h(v) = g^{-1}(v)$ the inverse of $g$ and set
$$
 H(\lambda, v) = \frac{\left(1-\tfrac{1}{\lambda}\right)v}{h(\lambda v)}. 
$$
Observe that $H(\lambda, v)$ is a continuous positive function  that converges to $0$ as $\lambda \to 1$ or $\lambda \to \infty$. 
If $h$ is everywhere differentiable  the maximum of $H(\lambda, v)$  for a fixed value of $v$ is therefore attained  for some $\lambda \in (1, \infty)$ satisfying 
$$
 \frac{\partial}{\partial \lambda } \left( \frac{(\lambda-1)v}{\lambda h(\lambda v)}\right) 
=
 \frac {h(\lambda v) + (\lambda-\lambda^2) v  h'(\lambda v)}{\left(\lambda h(\lambda v) \right) ^{2}} = 0,
$$
that is 
\begin{equation}\label{eqbestlabda}
 h(\lambda v) = \lambda (\lambda-1) v  h'(\lambda v).
\end{equation}
If this equation has a unique solution  $\lambda = \lambda(v)$, then \eqref{ineq2} holds for the  function 
$F(v) =  (1-\tfrac{1}{\lambda(v)})\frac{v}{h(\lambda(v)  v)}\ $.

\medskip

With this preparation, it is easy to prove the first statement in Corollary \ref{cortomain}:  Assume  that $\gamma_S(n) \geq g(n+1)$ with $g(r) = Cr^d$.
Then $h(v) = g^{-1}(v) = h(v) = (v/C)^{1/d}$ and $h'(v) = \frac{1}{vd}h(v)$. Equation \eqref{eqbestlabda} is then  the following equality:
$$
 0 =  h(\lambda v) + (\lambda-\lambda^2) v  h'(\lambda v) = h(\lambda v) \left( 1 + (\lambda-\lambda^2) \frac{1}{\lambda d}\right).
$$
The unique solution is $\lambda = d+1$, therefore \eqref{ineq2} holds with the function
$$
  F(v)  =  (1-\tfrac{1}{\lambda(v)})\frac{v}{h(\lambda(v)  v)} =  \frac{C^{\frac{1}{d}} d v^{\frac{d-1}{d}}}{(d+1)^{\frac{d+1}{d}}} ,
$$
which is equivalent to the inequality \eqref{isoppol}.

\medskip

We now prove the second statement  in Corollary \ref{cortomain}. Suppose  that $\gamma_S(n) \geq g(n+1)$   with $g(r) = C{{\rm e}^{b{r}^{\alpha}}}$,
then
\begin{equation}\label{hhexp}
 h(v) = g^{-1}(v)=  \left( {\frac {1}{b}\log  \left( {\frac {v}{C}} \right) } \right) ^{1/{\alpha}},
\end{equation}
and we have
$$
 h'(v) = \frac {1}{\alpha b^{1/\alpha} v}  \log \left( {\frac {v}{C}}  \right) ^{\frac{1}{\alpha}-1}
= \frac{h(v)}{\alpha v \log  \left(\frac {v}{C}\right)}.
$$
Equation  \eqref{eqbestlabda}   is then equivalent to 
\begin{equation}\label{eqbestlabdaexp1}
  \alpha \log \left(\frac{\lambda v}{C}\right) = (\lambda -1),  \quad  \text{equivalently}  \quad   v =  \frac{C}{\lambda}   \mathrm{e}^{\frac{\lambda-1}{\alpha}}.
 \end{equation}

We claim that for any $v \geq C \cdot  \dfrac{1}{\alpha} \mathrm{e}^{\frac{\alpha-1}{\alpha}} $, there exists a unique solution  $\lambda = \lambda(v)$ satisfying  \eqref{eqbestlabdaexp1}. It is indeed easy to check that the function 
$$
 f(x) =  C \cdot  \frac{1}{x} \mathrm{e}^{\frac{x-1}{\alpha}}
$$
has  its minimum value at $x = \alpha$ and  defines a diffeomorphism \ $ f : [\alpha, \infty) \to [ C \cdot  \frac{1}{\alpha} \mathrm{e}^{\frac{\alpha-1}{\alpha}} , \infty)$.  Clearly  $\lambda = f^{-1}(v)$ is the unique solution of \eqref{eqbestlabdaexp1}.
The function in \eqref{defF} is then given by 
$$
 F(v) =  \left(1 - \frac{1}{\lambda (v)}\right)  \frac{v }{h(\lambda(v) \cdot  v)}, 
$$
where $h(v)$ is defined in \eqref{hhexp}. 
This function can also be written as
\begin{eqnarray*}
 F(v) 
= \left(1 - \frac{1}{\lambda (v)}\right)  \frac{v}{ \left( \frac {1}{b} (\log (v) + \log(\lambda(v)) - \log(C) \right)^{1/\alpha}}
= \frac{ \mu(v)   v}{\left( \frac {1}{b} (\log (v)\right)^{1/\alpha}},
\end{eqnarray*}
with
\begin{equation}\label{defmu}
 \mu(v) =  \left(1 - \frac{1}{\lambda (v)}\right)    \left(1 + \frac{\log(\lambda (v))}{\log (v)}  - \frac{\log(C)}{\log(v)}\right)^{-1/\alpha}.
\end{equation}
Applying  Corollary \ref{corMainIsop}, we then have 
\begin{equation*} 
  \Card(\partial D)  \geq F(\Card(D))  = \mu(\Card(D))  \cdot   \frac{ \Card(D)}{ \left( \frac {1}{b} (\log (\Card(D))\right)^{1/\alpha}}.
\end{equation*}
It remains to  study the limit of $\mu(v)$ as $v \to \infty$. To see this, observe first that  equation \eqref{eqbestlabdaexp1} can also be written as 
\begin{equation}\label{eqbestlabdaexp2}
 \log(v)  = \frac{1}{\alpha}(\lambda-1)  -  \log(\lambda) + \log (C).
\end{equation}
From this equality, we see that $\lim_{v \to \infty}\lambda(v) = \infty$ and therefore.
$$
 \lim_{v \to \infty} \frac{\log(\lambda(v))}{\log{(v)}} =
 \lim_{\lambda  \to \infty} \frac{\log(\lambda)}{\left(\frac{1}{\alpha}(\lambda-1)  -  \log(\lambda) + \log (C)\right) } =  0.
$$
Therefore
\begin{equation}\label{limlambdav}
 \lim_{v \to \infty} \frac{\lambda(v)}{\log{(v)}}  =  {\alpha}  \quad \text{and} \quad  \lim_{v \to \infty} \frac{\log(\lambda(v))}{\log{(v)}} = 0.
\end{equation}
Applying  \eqref{limlambdav} to  \eqref{defmu}, we find that 
$
 \lim_{v\to \infty} \mu(v) = 1,
$
which completes the proof of the Corollary.
\qed

\section{Final Remarks}

\textbf{Remarks}  \textbf{1.}  The function $\lambda(v)$ in the proof of the second part of the Corollary is defined implicitly to be the 
solution of \eqref{eqbestlabdaexp1}.  We can also define  $\lambda(v)$ explicitly using the \emph{Lambert $W$-function}.
The Lambert function is multivaluate and in our case we consider the branch  $W_{-1} : [-\frac{1}{e}, 0)  \to (-\infty, -1]$ 
defined by the condition $W_{-1}(x) \exp(W_{-1}(x)) = x$, see  \cite{Beardon, CGHJK} for more on this function. 
It is then not difficult to check that 
$$
  \lambda(v) =  -\alpha W_{-1}\left(-\frac{C\mathrm{e}^{-\frac{1}{\alpha}} }{\alpha v}\right).
$$
Combining this formula with \eqref{defmu}, we see that $\mu(v)$ as the following explicit, albeit complex, expression
{\small 
\begin{equation}\label{defmu2}
 \mu(v) =  \left(1 + \frac{1}{\alpha W_{-1}\left(-\frac{C\mathrm{e}^{-\frac{1}{\alpha}} }{\alpha v}\right)}\right)    \left(1 + \frac{\log\left(-\dfrac{\alpha}{C} \cdot  W_{-1}\left(-\frac{C\mathrm{e}^{-\frac{1}{\alpha}} }{\alpha v}\right)\right) }{\log (v)}\right).
\end{equation}
}
\medskip

\textbf{2.} We have worked with the interior definition of the boundary of a subset $D\subset \Gamma$. We could have also chosen the 
\textit{exterior boundary}
$$
 \partial_S' D = \{x\in \Gamma \setminus D \mid \operatorname{dist}_S(x, D) = 1\}
=  \{x \in \Gamma \setminus  D \mid \exists s \in S \  \text{ such that } \ sx \in D\}.
$$
Because every point in $\partial_S D'$ is at distance $1$ from a point in $\partial_S D$ and vice-versa, we clearly have
$$
 \frac{1}{\Card(S)} \Card(\partial_S D) \leq \Card(\partial_S' D) \leq \Card(S) \Card(\partial_S D).
$$
The proof of  Theorem 1, and therefore all the estimates in the paper, also  holds if we replace the interior boundary $ \partial_SD$ with the exterior boundary  $\partial_S' D$. Only minor changes have to be made in step (i) of the proof; namely the function $f$ in \eqref{deff} has to be modified to
\begin{equation} \label{deff2}
 f(x) = y_mx, \quad \text{ where } m = \max \{ k \in \N \mid k \leq  n  \text{ and } y_k x \in \partial_S' D\}.
\end{equation}
In that case the inverse image of a point $z \in \partial_S' D$ satisfies
\begin{equation}\label{preimageZ2}
  f^{-1} (z) \subset \{y_1^{-1}z, \dots , y_{n}^{-1}z\},
\end{equation}
and we still have  $\Card (f^{-1}(z)) \leq n$ for any $z \in \partial_S' D$.

\medskip

\textbf{3.}  The function $\phi_S$ has a slightly different definition in \cite{PSC} and \cite{CSC}. These authors use instead the function 
\begin{equation}\label{def.phi2}
 \tilde{\phi}_S(t) = \min \{n \in \N \mid \gamma_S(n) > t\}.
\end{equation}
Observe that $\tilde{\phi}_S(t) \geq {\phi}_S(t)$, therefore the isoperimetric inequality \eqref{mainineq} still holds if one uses
$\tilde{\phi}_S$ instead of ${\phi}_S$, and it is in fact slightly weaker in that case.

\medskip

\textbf{4.} The isoperimetric inequality in Corollary   \ref{corMainIsop} is  not  optimal,  but  by construction, it is the best possible inequality that can be proved for a general finitely generated group based on the arguments we have used. It seems that proving sharper isoperimetric inequalities  for general groups would require new and possibly vers different ideas and techniques.

\rule{2cm}{0.1mm}%

\makeatletter 
{\small EPFL, Institut de Mathématiques, Station 8, 1015 Lausanne, Switzerland. \\ 
bruno.santoscorreia@epfl.ch, marc.troyanov@epfl.ch}
\makeatother %

\end{document}